%% file: moduli.tex
\documentclass[12pt, reqno]{amsart}
\usepackage{amssymb}
\usepackage{amsmath}
\usepackage[utf8]{inputenc}
\usepackage{longtable}
\usepackage{cite}
\usepackage[all,cmtip]{xy}
\usepackage{fullpage}

\include{pack}
	\title[The field of moduli]{The field of moduli of singular K3 surfaces}
	\author{Roberto Laface}
	\address{Institut f\"{u}r Algebraische Geometrie, Leibniz Universit\"{a}t Hannover, Welfengarten 1 30167 Hannover (Germany)}
	
	\email{laface@math.uni-hannover.de}
	
	\curraddr{Technische Universit\"at M\"unchen, Zentrum Mathematik M11 - Boltzmannstra{\ss}e 3, 85748 Garching bei M\"unchen (Germany)}
	
	\email{laface@ma.tum.de}
	
\begin{document}

\maketitle
\thispagestyle{empty}

\begin{abstract}
We study the field of moduli of singular K3 surfaces. We discuss both the field of moduli over the CM field and over $\Q$. We also discuss non-finiteness with respect to the degree of the field of moduli. Finally, we provide an explicit approach to the computation of the field of moduli.
\end{abstract}


\section{Introduction}

In algebraic geometry, it is sometimes convenient to work with general objects inside a moduli space, as the genericity assumption often seems to grant more control on the geometry. For instance, a general hypersurface of degree $d \geq 4$ in $\IP^3$ has Picard group of rank one, generated by the class of a hyperplane section. 

On the other extreme, one might want to consider objects with extra structure, and this often leads to very interesting arithmetic applications. As an example, the generic elliptic curve will have endomorphism ring isomorphic to $\Z$, but it is well-known that elliptic curves with complex multiplication (CM), meaning those whose endomorphism ring is strictly bigger than the ring of integers, have attracted a great deal of researchers in algebraic and arithmetic geometry. In the moduli space of (smooth) elliptic curves, CM elliptic curves are a countable set, and it can be shown that they correspond to quadratic imaginary numbers (see, for instance, \cite[Ch.~IV, Sec.~4]{hartshorne77}).\\

Abelian surfaces and K3 surfaces, which are the natural generalization of elliptic curves to dimension two, are no exception: the groundbreaking work of Shioda and Mitani \cite{shioda-mitani74} has revealed a deep connection between the geometry of \textit{singular} abelian surfaces (abelian surfaces with maximum Picard number) and the arithmetic of quadratic forms. Also, this connection extends to singular K3 surfaces, as it was shown by Shioda and Inose \cite{shioda-inose77}, by means of what Morrison later called Shioda-Inose structures \cite{morrison84}. 

The arithmetic data of a singular abelian surface is encoded in its transcendental lattice, and a Shioda-Inose structure associates to it a singular K3 surface with the same transcendental lattice, thus preserving the arithmetic information. This has been employed, for instance, by Sch\"utt in the study of the field of definition of singular K3 surfaces \cite{schuett07}: he proved that a singular K3 surface $X$ always admits a model over a ring class field $H/K$, $K$ being the field $K=\Q(\disc\rT(X))$, generalizing previous results of Inose \cite{inose78}. Morever, generalizing previous work of Shimada \cite{shimada09}, he describes the conjugate varieties of $X$ (modulo $\C$-isomorphism) under the action of $\Aut(\C/K)$: this is done by looking at the corresponding transcendental lattices, and it is best understood in the language of genus theory of quadratic forms.

This suggests that, given a good notion of field of moduli, the degree of the field of moduli should be exactly the number of Galois conjugates of $X$. However, apart from the aforementioned result, nothing is known in general for the field of moduli of singular K3 surfaces. This paper aims at describing the field of moduli of singular K3 surfaces.\\

By using an idea of \v{S}afarevi\v{c} \cite{shafarevich96}, we reduce the problem of studying the field of moduli of a singular K3 surface $X$ to the study of the analogous field of a singular abelian surface $A$ with transcendental lattice $\rT(A)=\rT(X)$ (this condition can always be achieved by means of Shioda-Inose structure). Since the transcendental lattice of $X$ and $A$ are isomorphic, they will have the same field of moduli. Our main tools are Galois theory, the theory of complex multiplication on elliptic curves and the theory of quadratic forms. 

Throughout the paper, we stress the analogy between CM elliptic curves and singular K3/abelian surfaces. For instance, singular K3/abelian surfaces are always defined over a number field (see, for example, \cite{serre-tate68} and \cite{serre-tate68}). It turns out that, in the case of K3/abelian surfaces, we recover a very similar picture to the case of CM elliptic curves and their field of moduli.\\

\noindent
We give a short overview of the results achieved in the present article. We start by giving a notion of \textit{field of $K$-moduli}, where $K$ is an arbitrary field. Typically, we will consider the cases where $K$ is the CM field of our singular K3 surface and the case $K=\Q$.

We first show that, when $K$ is the CM field of a singular K3 surface $X$, the field of $K$-moduli $M_K$ is a Galois extension of $K$ of degree $g$, $g$ being the order of the genus of the transcendental lattice of $X$ seen as a quadratic form.

Then, we study the field of $\Q$-moduli of $X$, say $M_\Q$: it turns out that it is a degree $g$ extension of $\Q$, having index two inside $M_K$. In general, $M_\Q/\Q$ is not a Galois extension, as we point out in the examples: it is enough to consider a singular K3 surface with class number three.

Afterwards, we investigate non-finiteness of singular K3 surfaces with respect to the field of moduli, and how the field of $K$-moduli changes as we vary the singular K3 surface: as an example, we prove that the field of moduli of a singular K3 surface is independent of the index of primitivity of the transcendental lattice. We also provide an explicit description of the field of moduli that can be implemented on a computer algebra system. 

\subsection*{Acknowledgements}
It is a pleasure to thank Matthias Sch\"utt for fruitful conversations and for sharing his insights on this subject.

\section{Preliminaries}

\subsection{Singular surfaces}
We start by recalling some basics on singular surfaces, and, to this end, let us work over the field $\C$ of complex numbers. If $X$ is a smooth algebraic surface, we can define the N\'{e}ron-Severi lattice of $X$:
it is the group of divisors on $X$, modulo algebraic equivalence, namely
\[\NS(X) := \Div(X) / \sim_{\rm alg},\]
together with the restriction of the intersection form on $\rH^2(X,\Z)$. Its rank $\rho(X):=\rank \NS(X)$ is called \textit{Picard number} of $X$; the Picard number measures how many different
curves lie on a surface. By the Lefschetz theorem on $(1,1)$-classes, we have the bound
\[\rho(X) \leq h^{1,1}(X) = b_2(X) - 2p_g(X),\]
where $b_2(X) := \rank \rH^2(X,\Z)$ and $p_g(X) : =\dim_\C \rH^0(X,\omega_X)$. 

We can consider the lattice 
\[ \rH^2(X,\Z)_\text{free} := \rH^2(X,\Z)/(\text{torsion}),\]
and since $\NS(X) \subset \rH^2(X,\Z)$, also $\NS(X)_\text{free} \subset \rH^2(X,\Z)_\text{free}$; $\NS(X)_\text{free}$ is a lattice of signature $(1,\rho(X)-1)$. Its orthogonal complement
$\rT(X) \subset \rH^2(X,\Z)_\text{free}$ is called the \textit{transcendental lattice} of $X$, and it has signature 
\[(2p_g(X), h^{1,1}(X) - \rho(X)).\]

A smooth algebraic surface with maximum Picard number, i.e.~$\rho(X) = h^{1,1}(X)$, is called a \emph{singular surface}. In this case, the
transcendental lattice acquires the structure of a positive definite lattice of rank $2p_g(X)$. Both in the case of singular abelian surfaces and singular K3 surfaces, $\rT(X)$ is a positive definite rank-two lattice. The exponential sequence
\[0 \longra \Z \longra \ko_A \longra \ko_A^\times \longra 0\]
yields a long exact sequence in cohomology, from which we can extract a map
\[p_A \, : \, \rH^2(A,\Z) \longra \rH^2(A,\ko_A) \cong \C,\]
since $p_g(A)=1$; the map $p_A$ is called \textit{the period of $A$}. For further details on singular surfaces, see \cite{shioda-mitani74} and \cite{beauville14}.

\subsection{Class group theory}
We recall a few facts on integral binary quadratic forms; for a detailed account, the reader is suggested to see \cite{cox13}. 
Given a form 
\[Q(x,y)=ax^2+bxy+cy^2\]
the quantity $\gcd(a,b,c)$ is called \textit{index of primitivity} of $Q$, and $Q$ is said \emph{primitive} if $\gcd(a,b,c)=1$. Sometimes, it is convenient to extract the \textit{primitive part} of a form $Q$: this is the quadratic form $Q_0$ such that $mQ_0 =Q$, $m$ being the index of primitivity of $Q$. 

A form $Q$ \emph{represents} $m \in \Z$ if $m=Q(x,y)$ for some
$x,y \in \Z$; if moreover $\gcd(x,y)=1$, then we say that $Q$ \emph{properly represents} $m \in \Z$. A quadratic form $Q$ as above
will be denoted in short by $Q=(a,b,c)$. Two forms $Q = (a,b,c)$ and $Q'=(a',b',c')$ are equivalent (properly equivalent, respectively) if there exists
$\begin{pmatrix} p & q \\ r & s\end{pmatrix} \in \GL_2(\Z)$ ($\SL_2(\Z)$, respectively) such that 
\[Q(px+qy,rx+sy) = Q'(x,y).\]

The \textit{discriminant} of a form $Q=(a,b,c)$ is the integer $D:=b^2-4ac$. The set of proper equivalence classes of primitive forms of discriminant $D$ is called the \emph{(form) class group of discriminant $D$},
and it is denoted by $C(D)$; we will denote the class of a form $Q$ by $[Q]$. The class group is equipped with the Dirichlet composition
of forms: by \cite[Lemma 3.2]{cox13}, if $Q=(a,b,c)$ and $Q'=(a',b',c')$ are primitive forms of discriminant $D$, such that
\[\gcd\Big(a,a',\frac{b+b'}{2}\Big)=1,\]
then the composition $Q*Q'$ is the form $(aa',B,C)$, where $C=\frac{B^2-D}{4aa'}$ and $B$ is the integer,
unique modulo $2aa'$, such that
\[\left.
\begin{cases}
B \equiv b \mod 2a, \\
B \equiv b' \mod 2a', \\
B^2 \equiv D \mod 4aa'.
\end{cases}\right.\]
Naturally, we put $[Q]*[Q']:=[Q*Q']$. 

Recall that, fixed a quadratic imaginary field $K$, an order $\ko$ is a subring of $K$ containing the unity of $K$ which has also the 
structure of a rank-2 free $\Z$-module; every order $\ko$ can be written in a unique way as
\[\ko=\Z + fw_K \Z, \quad w_K:=\frac{d_K+\sqrt{d_K}}{2}, \quad d_K:=\disc \ko_K, \quad f \in \Z^+.\]
The integer $f$ is called the \emph{conductor} of $\ko$, and it characterizes $\ko$ in a unique way; we will denote the order of conductor $f$ in $\ko_K$ by $\ko_{K,f}$.

For an order $\ko$ in a quadratic field $K$, it is possible to define a class group $C(\ko)$: letting $I(\ko)$ denote the group of proper
fractional ideals, meaning those whose CM ring is $\ko$ itself, and letting $P(\ko)$ be the subgroup generated by the principal ones, we
set $C(\ko) := I(\ko)/P(\ko)$, and we call it the \emph{ideal class group} of $\ko$. An important result in algebraic number theory states
that if $\disc \ko = D$, then $C(D) \cong C(\ko)$; from now on, we will use interchangeably the two class groups to our convenience. The order of the class group $C(\ko)$ is called the \emph{class number} of $\ko$, and it is denoted by $h(\ko_{K,f})$.

\subsection{Two interesting spaces of singular surfaces}
Let $\Sigma^{\rm Ab}$ be the set of moduli of singular abelian surfaces, i.e.~isomorphism classes of singular abelian surfaces; in \cite{shioda-mitani74}, Shioda and Mitani described $\Sigma^{\rm Ab}$ by means of the transcendental lattice $\rT(A)$ associated to any singular abelian surface $A$. We say that an ordered basis $\lbrace t_1, t_2 \rbrace$ of T(A) is positive if
\[\text{Im} ( p_A(t_1) / p_A(t_2) ) > 0,\]
and T(A) with a choice of a positive basis is said to be positively oriented.

Notice that the transcendental lattice $\rT(A)$ is an even lattice, and after choosing a basis one has
that
\[ \rT(A) = \begin{pmatrix} 2a & b \\ b & 2c \end{pmatrix}, \qquad a,c>0, \qquad b^2-4ac. \]
Thus we can always associate to it the quadratic form $(a,b,c)$. This realizes a 1:1 correspondence, and therefore we can naturally see the transcendental lattice as a integral binary quadratic form. 

We can associate to any quadratic form $Q=(a,b,c)$ an abelian surface $A_Q$. In order to describe the correspondence, we set
\[\tau(Q):=\dfrac{-b+\sqrt{D}}{2a}, \qquad D:= \disc Q = b^2-4ac,\]
and we will denote by $E_{\tau}$ the elliptic curve $\C / \Lambda_{\tau}$, $\Lambda_{\tau}$ being the lattice $\Z + \tau\Z$.
The abelian surface associated to a form $Q$ is then defined as the product surface
\[A_Q:=E_{\tau} \times E_{a\tau+b},\]
where $\tau = \tau(Q)$. 

The mapping $Q \mapsto A_Q$ realizes a 1:1 correspondence between $\SL_2(\Z)$-conjugacy classes of binary forms and isomorphism classes
of singular abelian surfaces, namely
\[ \Sigma^{\rm Ab} \longlra \kq^+ / \SL_2(\Z),\]
$\kq^+$ being the set of positive definite integral binary quadratic forms. By forgetting the orientation, we get a 2:1 map
$\Sigma^{\rm Ab} \longra \kq^+/\GL_2(\Z)$, which is just taking the transcendental lattice of an abelian surface:
\[\Sigma^{\rm Ab} \ni [A] \longmapsto [\rT(A)] \in \GL_2(\Z).\]
As a consequence, we get that every singular abelian surface $A$ is isomorphic to the product of two isogenous elliptic curves with
complex multiplication.

One can build the analogous space $\Sigma^{\rm K3}$ of singular K3 surfaces, and ask for its structure. In their paper, Shioda and Mitani \cite{shioda-mitani74} showed that by taking the Kummer surface of a singular abelian surface, one is able to recover all singular K3 surfaces whose transcendental lattice has primitivity index which is divisible by 2. Later, Shioda and Inose \cite{shioda-inose77} proved the surjectivity of the period map for singular K3 surfaces by means of Shioda-Inose structures (as Morrison called them in \cite{morrison84}): if $A$ is an abelian surface, a Shioda-Inose structure associated to $A$ a K3 surface $X = SI(A)$, which is a 2:1 cover of $\Km (A)$ and has the property that $\rT(X) = \rT(A)$.
\[A \xdashrightarrow{ \quad 2:1 \quad} \Km(A) \xdashleftarrow{\quad 2:1 \quad} X \]

It turns out that also (isomorphism classes of) singular K3 surfaces are uniquely characterized by their transcendental lattice, and thus $\Sigma^{\rm Ab} \cong \Sigma^{\rm K3}$; in particular, this implies that two singular abelian surfaces are isomorphic if and only if the corresponding K3 surfaces via a Shioda-Inose structure are. This amounts to saying that we can interchangeably consider singular abelian surfaces and singular K3 surfaces when studying the field of moduli.

\subsection{CM theory of elliptic curves}
We recall a couple of elementary facts about the CM theory of elliptic curves; for a reference, see \cite[Ch.~2]{silverman09}. Let $\ke ll (\ko)$ be the set of isomorphism classes of elliptic curves with CM by the order $\ko \subset K$. Since a proper $\ko$-ideal is also a lattice, quotienting by $\ko$-ideals induces a map
\[C(\ko) \longra \ke ll(\ko), \qquad \bar{\ga} \longmapsto [\C / \ga],\]
which is an isomorphism. Multiplication of ideal classes and lattices gives a simply transitive action
\[C(\ko) \times \ke ll(\ko) \longra \ke ll(\ko), \qquad (\bar{\ga},[\C / \Lambda]) \longmapsto \bar{\ga} \ast [\C/\Lambda]:=[\C/\ga^{-1} \Lambda].\]
Another action on $\ke ll(\ko)$ is given by the absolute Galois group $\Gal(\bar{K}/K)$:
\[\Gal(\bar{K}/K) \times \ke ll(\ko) \longra \ke ll(\ko), \qquad (\sigma,[E]) \longmapsto [E^\sigma].\]

Now, let us fix $[E] \in \ke ll(\ko)$; given $\sigma \in \Gal(\bar{K}/K)$, we can form $[E^\sigma]$, and, by using the action of $C(\ko)$, there exists a unique $\bar{\ga} \in C(\ko)$ such that $\ga \ast [E] = [E^\sigma]$. This correspondence defines a surjective homomorphism
\[F: \Gal(\bar{K}/K) \longra C(\ko), \qquad \sigma \longmapsto F(\sigma) \, : \, F(\sigma) \ast [E] = [E^\sigma].\]
One of the properties of this map is that it is independent of the curve $[E]$ chosen to define it, and thus we have that $F(\sigma) \ast [E] = [E^\sigma]$, $\forall \sigma \in \Gal(\bar{K}/K)$ and $\forall [E] \in \ke ll(\ko)$. 

Notice that the action 
\[C(\ko) \times \ke ll(\ko) \longra \ke ll(\ko)\]
can be interpreted in terms of quadratic forms. Indeed, to any $[E] \in \ke ll(\ko)$, one can associate a quadratic form $Q$ such that $j(\tau(Q)) = j(E)$. Then, the action is isomorphic to the action 
\[C(\ko) \times C(\ko) \longra C(\ko), \qquad (\bar{\ga},\bar{\gb}) \longmapsto \bar{\ga}^{-1} \bar{\gb}.\]
Also, by class group theory, we can phrase everything in terms of the corresponding classes of quadratic forms, where now multiplication of ideal classes corresponds to the Dirichlet composition. Under this interpretation, the map 
\[F : \Gal(\bar{K}/K) \longra C(\ko), \qquad \sigma \longmapsto F(\sigma),\]
where $F(\sigma)$ is the element of $C(\ko)$ such that
\[[E^\sigma] = F(\sigma) \ast [E] = F(\sigma) \ast [\C/\ga] = [\C / F(\sigma)^{-1}  \ga] = [F(\sigma)]^{-1} \ast [Q],\]
where $[Q]$ corresponds to $[E]$, and $[F(\sigma)]$ is meant as the class of quadratic forms corresponding to $F(\sigma) \in C(\ko)$.

\subsection{The group of id\'eles}
Given a number field $K$, we denote by $I_K$ the group of fractional ideals in $K$. We can define the so-called \textit{group of id\'{e}les} by setting 
\[ \II_K := \Big\lbrace (a_v) \in \prod_v K_v^\times \, \Big\vert \, a_v \in \ko_v^\times \text{for all but finitely many $v$}\Big\rbrace.\]
There is a canonical surjective homomorphism $\id$
\[\II_K \longra I_K, \qquad (a_v) \longmapsto \prod_{\text{$v$ finite}}p_v^{\ord_{p_v}(a_v)},\]
which associates to every id\'ele an element of $I_K$. There is also a canonical injective (diagonal) homomorphism
\[K^\times \longra \II_K, \qquad a \longmapsto (a,a,a,\dots ),\]
with discrete image. 

The statement of the main theorems of class field theory in terms of ideals is very explicit. However, it has the big disadvantage of working for a fixed modulus $\gm$ at the time, and so it describes only the abelian extensions whose conductor divides $\gm$. On the other hand, the statements in terms of id\'{e}les allow one to consider infinite abelian extensions, or equivalently all finite abelian extensions simultaneously. It also relates local and global class field theory, namely the global Artin map to its local components.

\begin{prop}
There exists a unique continuous surjective homomorphism $\phi_K : \II_K \longra \Gal(K^{\rm ab}/K)$ with the following property: for any $L \subset K^{\rm ab}$ finite over $K$ and any prime $w$ of $L$ lying over a prime $v$ of $K$, the diagram
\[\xymatrix{
&K_v^\times \ar[d] \ar[rr]^{\phi_v} & &\Gal(L_w/K_v) \ar[d] \\
&\II_K \ar[rr]^{\phi_{L/K}} & &\Gal(L/K)
}
\]
where the bottom map sends $a \in \II_K$ to $\phi_K(a)\vert_L$. Here, $\phi_v$ is the local component of the Artin map at the place $v$.
\end{prop}
In particular, for any finite extension $L$ of $K$ which is contained in $K^{\rm ab}$, $\phi_K$ gives rise to a commutative diagram
\[\xymatrix{
&\II_K \ar[rr]^{\phi_{K}} \ar[rrd]_{\phi_{L/K}} & &\Gal(K^{\rm ab}/K) \ar[d]^{\vert_L}\\
& & &\Gal(L/K)
}
\]
1
We recall that the Main Theorem of Complex Multiplication makes use of the group of id\'{e}les $\II_K$ to control the Galois conjugates of an elliptic curve with CM in $K$. Let $K$ be an imaginary quadratic field and $E$ an elliptic curve with CM in $K$; then, there exist an order $\ko \subset K$ and a fractional ideal $\ga \subset \ko$ such that $E \cong \C / \ga$, and thus $E$ has CM in the order $\ko$.

\begin{thm}[Main Theorem of Complex Multiplication, Theorem 5.4 of \cite{shimura71}]
Let $E = \C/\Lambda$ be an elliptic curve with CM by an order in $K$. Let $\sigma \in \Gal(\bar{K}/K)$ and $s \in \II_K$ such that $\sigma = \phi_K(s)$ on $K^{\rm ab}$. Then, there exists an isomorphism
\[E^\sigma \cong \C / s^{-1} \Lambda.\]
\end{thm}

\section{The field of $K$-moduli}\label{section_3}

\subsection{A new definition}
We define the \textit{field of $K$-moduli} $M_K$ of a variety $X$, where $K$ is a given field. This field was first introduced by Matsusaka \cite{matsusaka56} as the \textit{relative field of moduli} (or \textit{field of moduli} over $K$), and it was defined to be the intersection of all fields of definition of $X$ which contain $K$, in other words
\[M_K:=\bigcap_{\substack{\text{$X$ defined over $L$} \\ L \supset K}} L.\]

Later, Koizumi \cite{koizumi72} adjusted the definition to positive characteristic geometry by adding the extra condition that for an automorphism $\sigma \in \Aut(\Omega / K)$, where $\Omega$ is a fixed universal domain\footnote{Given a field $K$, a universal domain $\Omega$ is an extension of $K$ with infinite transcendence degree over $K$. Universal domains were the fundamental object algebraic geometry was based on before the advent of Grothendieck. More details can be found in the fundational book of Weil \cite{weil46}; this uses notions very much different from the modern language of schemes and it is quite hard to read at times.},
\[\sigma \in G:=\lbrace \sigma \in \Aut(\Omega/K) \, \vert \, X^\sigma \in [X] \rbrace \Longleftrightarrow \sigma_{\vert M_K} = \id_{M_K},\]
where by $[X]$ we denote the isomorphism class of $X$. For our purposes, it is best to introduce the following

\begin{defi}
	The \textit{field of $K$-moduli} of $X$ is the subfield of $\C$ fixed by the group
	\[G:=\lbrace \sigma \in \Aut(\C/K) \, \vert \, X^\sigma \in [X] \rbrace.\]
\end{defi}

In practice, we are dropping Matsusaka's condition and keeping the one Koizumi introduced. Notice that, unlike in the case of Koizumi's definition \cite{koizumi72}, our field of moduli always exists and it is unique by Galois theory. Following \cite{koizumi72}, if the characteristic of the ground field is zero, then $M_K$ is contained in any field of definition for $X$ which contains $K$, and thus we have the following extension
\[M_K \subset \bigcap_{\substack{\text{$X$ defined over $L$} \\ L \supset K}} L,\]
which in fact is algebraic and Galois. We remark that the right-hand side of this inclusion is quite a mysterious object in general.

If $X$ is a variety, by the \textit{absolute field of moduli} of $X$ we will mean the field of $\Q$-moduli, i.e.~the field $M_\Q$ such that for all automorphisms $\sigma \in \Aut(\C/\Q)$,
\[X^\sigma \in [X] \Longleftrightarrow \text{$\sigma$ acts trivially on $M_\Q$};\]
equivalently, it is defined as the fixed field of the group
\[G:=\lbrace \sigma \in \Aut(\C/\Q) \, \vert \, X^\sigma \in [X] \rbrace.\]
Galois theory once again guarantees that this field is unique for a given variety $X$.

If $X$ is a variety and $\tau \in \Gal(\C / \Q)$, let $X^\tau$ denote the variety obtained by conjugating $X$ by $\tau$. Suppose we want to study the field of $L$-moduli, for some number field $L$, and denote by $G(X)$ (respectively, $G(X^\tau)$) the group fixing the modulus of $X$ (respectively, $X^\tau$) and by $M(X)$ (respectively, $M(X^\tau)$) the field of $L$-moduli. Then, one can show that:
\begin{enumerate}
	\item $G(X)$ only depends on the isomorphism class of $X$;
	\item $G(X^\tau) = \tau \cdot G(X) \cdot \tau^{-1}$;
	\item $M(X^\tau) = \tau ( M(X))$.
\end{enumerate} 

\subsection{A little motivation}
Let $X$ be a singular K3 surface, and let
\[\rT(X) \cong \begin{pmatrix} 2a & b \\ b & 2c \end{pmatrix} = (a,b,c)\]
be its transcendental lattice, where the right-hand side equality identifies $\rT(X)$ with the corresponding quadratic form. To $\rT(X)$, one can associate two gadgets: the first one is $\det \rT(X)$, and the second one is $\disc \rT(X)$, when $\rT(X)$ is regarded as a quadratic form. Clearly,
\[\det \rT(X) = - \disc \rT(X) < 0,\]
and thus $K := \Q(\sqrt{\det \rT(X)})$ is a quadratic imaginary field. We will call $K$ the \textit{CM field}\footnote{There are other notions of CM field currently in use. For example, a number field $K$ is a CM field if it is a totally imaginary quadratic extension of a totally real field. In fact, the CM field of a singular K3 surface is also a CM field in the latter sense.} of $X$. If $X$ is the K3 surface associated to a singular abelian surface $A$ via a Shioda-Inose structure (so that, in particular, $\rT(A) = \rT(X)$), we will say that $K$ is the CM field of $A$ as well. Generalizing a previous result of Shimada \cite{shimada09}, Sch\"{u}tt was able to prove the following result

\begin{thm}[Theorem 5.2 in \cite{schuett07}]\label{genus}
	Let $X$ be a singular K3 surface, and let $\rT(X)$ be its transcendental lattice. Assume that $X$ is defined over a Galois extension $L/K$. Then, the action of the Galois group $\Gal(L/K)$ spans the genus of $\rT(X)$, i.e.
	
	\[\big({\rm genus\ of\ } \rT(X) \big) = \big\lbrace [\rT(X^\sigma)] \, : \, \sigma \in \Gal(L/K) \big\rbrace.\]
\end{thm}

Here, as the genus is defined for primitive quadratic forms only, we mean the following: consider the primitive part of $\rT(X)$, so that $\rT(X) = m \rT(X)_0$, $m$ being the index of primitivity of $\rT(X)$. Then, 
\[\big(\text{genus of }\rT(X) \big) = \big\lbrace m [T] \, : \, [T] \text{ lies in the genus of }\rT(X)_0 \big\rbrace.\]

Set $L:=H(\disc \rT(X))$ in Theorem \ref{genus}, where $H(D)$ denotes the ring class field of the order in $K$ of discriminant $D$, for $D < 0$. Galois theory tells tells us that
\[\Gal(L/\Q) \cong \Gal(L/K) \rtimes \Gal(K/\Q),\]
where $\Gal(K/\Q)$ accounts for the complex conjugation (for a reference, see \cite[Ch.~9]{cox13}). But complex conjugation has the effect of sending a singular K3 surface of transcendental lattice $(a,b,c)$ to the singular K3 surface with transcendental lattice $(a,-b,c)$, so it acts as inversion on the corresponding class group (see \cite{shioda-mitani74} and \cite{schuett07}). By observing that a form and its inverse lie in the same genus, we conclude that
\begin{align*}
(\text{genus of $T(X)$}) =& \lbrace [T(X^\sigma)] \, \vert \, \sigma \in \Gal(L/K) \rbrace =\\
=& \lbrace [T(X^\sigma)] \, \vert \, \sigma \in \Gal(L/\Q) \rbrace.
\end{align*}
This observation suggests a connection between the field of moduli of a singular K3 surface and the genus of its transcendental lattice, even in the case of the field of $\Q$-moduli.

The classification of decompositions of a singular abelian surface \cite{laface15} allows us to tell something more about the field of moduli of $X$ containing $K$. Recall that $M_K$ is contained in the intersection of all possible fields of definition for $X$. Then, by means of Shioda-Inose structures, we can study $X$ by means of those abelian surfaces $A$ whose transcendental lattice equals $\rT(X)$. Let $A$ be such a surface, and consider all product surfaces $E_1 \times E_2$ isomorphic to $A$ (which we know explicitly by \cite{laface15}); if $j_k := j(E_k)$, by work of Sch\"{u}tt \cite{schuett07}, $X$ admits a model over $\Q(j_1 j_2, j_1 +j_2)$. Therefore, considering all admissible pairs $(E_1,E_2)$ as above, we see that
\[M_K \subseteq \bigcap_{\text{$X$ defined over $L$}} L \subseteq \bigcap_{\text{$j_1,j_2$ as above}} \Q(j_1 j_2, j_1 +j_2).\]
We deduce a slightly clearer picture of what $M_K$ looks like, as we know where it has to sit as an extension of $\Q$. Namely, $M_K$ lies in right-hand side above, which is theoretically clear. In practice, describing it is a hard task, as this involves the computation of $j$-invariants.

\subsection{The case of elliptic curves}

Our toy example is the case of an elliptic curve $E$, for which one always has a Weierstra\ss\ model
\[ y^2 = x^3 + Ax + B,\]
for some $A,B \in \C$. It can be proven (see \cite[Ch.\,1]{silverman94}) that an elliptic curve $E$ can be defined over the field $\Q(j_E)$; moreover, the field of $\Q$-moduli of $E$ is again $\Q(j_E)$. 

Let now $E$ be a CM elliptic curve. The theory of complex multiplication tells us (see \cite[Ch.\,2]{silverman94}) that $j_E \in \overline{\Q}$, i.e.\,the $j$-invariant of a CM elliptic curve is always an algebraic number. Suppose that $E$ has CM by an order $\ko$ in $K = \Q(\sqrt{D})$. Then, by means of class field theory, one can show that there exists a commutative diagram of field extensions,

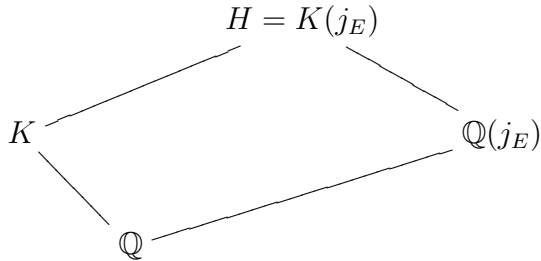
\[
\xymatrix{
	& & &H = K(j_E) & \\
	&K  \ar@{-}[rru] & & &\Q(j_E) \ar@{-}[lu] \\
	& &\Q \ar@{-}[lu] \ar@{-}[rru] & &
}
\]
where $H := H(D)$ is the ring class field corresponding to the order $\ko$ (for details, consult \cite{shimura71}). 

We would like to let the reader notice that $K(j_E)$ is indeed the field of $K$-moduli of $E$. Our study of the field of moduli in the rest of the chapter will reveal that this very picture carries over to singular K3 surfaces (and singular abelian surfaces). 

\subsection{An alternative definition of $M_K$}
As a singular K3 surface is defined over a number field by a result of Inose \cite{inose78}, when studying the field of moduli one would like to consider the field
\[ \bar{K}^{G'}, \qquad G':=\lbrace \sigma \in \Aut(\bar{K}/K) \, \vert \, X^\sigma \in [X] \rbrace,\]
rather than $\C^G$, as we defined it above. In fact, one has that $\C^G = \bar{K}^{G'}$; also this is independent of the fact that we are working on a singular K3 surface, as the following more general result shows.

\begin{prop}\label{FM}
	Let $X$ be a variety defined over a number field containing $K$. Then the fields $\C^G$ and $\bar{K}^{G'}$ coincide.
\end{prop}

\begin{proof}
	If $X$ is defined over a number field containing $K$, then $\Aut(\C/\bar{K})\subseteq G$, and thus $\C^G \subseteq \C^{\Aut(\C/\bar{K})}=\bar{K}$ (see \cite[Theorem 9.29]{milneFT}). This immedialtely implies that $\C^G \subseteq \bar{K}^{G'}$. The reverse inclusion follows from the surjectivity of the restriction map $\vert_K: G \longra G'$.
\end{proof}

As every singular K3 surface can be defined over a number field, we can define the field of $K$-moduli of a singular K3 surface to be the field 
\[M_K := \bar{K}^{G_K}, \qquad G_K:=\lbrace \sigma \in \Gal(\bar{K}/K) \, \vert \, X^\sigma \in [X] \rbrace.\]
In the following, we will be concerned with finding explicitly the group $G_K$, as it characterizes uniquely, thanks to Galois theory, the field of moduli.

\section{Characterization  in the primitive case}\label{prim}

\subsection{Statement of the result}
Let $X$ be a singular K3 surface, with transcendental lattice $\rT(X) = Q =mQ_0$ ($Q_0$ being the primitive part of $\rT(X)$), and discriminant $\disc \rT(X) = D = m^2 D_0$ ($D_0$ being the discriminant of $Q_0$). Recall that we can always find a singular abelian surface $A$ such that $X$ is obtained from $A$ by means of the Shioda-Inose structure, and in particular such that $\rT(A) = \rT(X)$. In light of this, notice that determining the field of moduli of $X$ is equivalent to determining the field of moduli of any such $A$, so that we can reduce to considering the problem for singular abelian surfaces.

We will now proceed in giving a different characterization of $G_K$. In what follows, let us assume additionally that $m=1$, which is to say that the transcendental lattice $\rT(X)$ is primitive. Under this assumption, for any decomposition $A \cong E_1 \times E_2$, the quadratic forms $Q_1$ and $Q_2$ corresponding to the elliptic curves $E_1$ and $E_2$ both lie in $C(D)\cong C(\ko)$, $\ko$ being the order of discriminant $D$. Observe that, if we fix a decomposition of $A\cong E_1 \times E_2$, then 
\[X^\sigma \in [X] \Longleftrightarrow A^\sigma \in [A] \Longleftrightarrow E_1^{\sigma} \times E_2^{\sigma} \cong E_1 \times E_2.\]
We will prove the following

\begin{thm}\label{thmprim}
	Let $X$ be a singular K3 surface with primitive transcendental lattice, and let $H$ be the ring class field of $\ko$, the order of discriminant $\disc \rT(X)$. Then the field of $K$-moduli is
	\[M_K = \bar{K}^{G_K}, \qquad G_K = (\vert_H)^{-1}\Gal(H/K)[2];\]
	it is a Galois extension of $K$ of degree $g$, $g$ being the order of the genus of the transcendental lattice of $X$.
\end{thm}

The proof is divided into two steps. First, we will prove that $G_K$ restricts to the subgroup of $2$-torsion elements of $\Gal(H/K)$, and thus it is a closed and normal subgroup of $\Gal(\bar{K}/K)$ with respect to the Krull topology. Afterwards, we will use these facts to study the field extension $M_K/K$, hence to prove Theorem \ref{thmprim}.

\subsection{The group $G_K$}
By the previous discussions, it follows that
\begin{align*}
G_K&=\lbrace \sigma \in \Gal(\bar{K}/K) \, \vert \, X^\sigma \in [X] \rbrace \\
&=\lbrace \sigma \in \Gal(\bar{K}/K) \, \vert \, E_1^{\sigma} \times E_2^{\sigma} \cong E_1 \times E_2 \rbrace.
\end{align*}
We will now proceed in giving a different characterization of $G_K$.

\begin{prop}\label{fundcond}
	$G_K = F^{-1}(C(\ko)[2])$.
\end{prop}
\begin{proof}
	Let $Q_i$ be the form corresponding to $E_i$ ($i=1,2$), and let $Q_i^\sigma$ be the one corresponding to $E_i^\sigma$ ($i=1,2$). By use of the map 
	\[F : \Gal(\bar{K}/K) \longra C(\ko),\]
	we get that
	\[ [Q_1^\sigma]=[F(\sigma)] ^{-1}  \ast [Q_1] \qquad \text{and} \qquad [Q_2^\sigma]=[F(\sigma)] ^{-1}  \ast [Q_2],\]
	where here we make use of the fact that $F$ is independent of the elliptic curve (and thus of the quadratic form) we use to define it. By \cite{laface15}, we see that
	\begin{align*}
	E_1^\sigma \times E_2^\sigma \cong E_1 \times E_2 &\Longleftrightarrow Q_1^\sigma \ast Q_2^\sigma = Q_1 \ast Q_2\\ &\Longleftrightarrow F(\sigma)^2 =1.
	\end{align*}
\end{proof}

There is a commutative diagram
\begin{equation}
\tag{$\dagger$}
\xymatrix{
	&\Gal(\bar{K}/K) \ar@{->>}[r]^{F} \ar@{->>}[d]_{{\vert_{H}}} & C(\ko) \\
	&\Gal(H/K) \ar[ru]_{\cong} &
}
\end{equation}
where $H:=H(\ko)$, which follows from class group theory and says that $F$ is an isomorphism on the restriction of the elements of $\Gal(\bar{K}/K)$ to $H$. In particular, $C(\ko)[2] \cong \Gal(H/K)[2]$, and thus 
\[ G_K = \lbrace \sigma \in \Gal(\bar{K}/K) \, : \, (\sigma\vert_H)^2 = \operatorname{id}_{H} \rbrace.\]
This implies the following

\begin{cor}
	$G_K = (\vert_H)^{-1}(\Gal(H/K)[2])$.
\end{cor}

We now turn to describing $G_K$ as a topological subgroup of $\Gal(\bar{K}/K)$.
\begin{prop}
	$G_K$ is a closed normal subgroup of $\Gal(\bar{K}/K)$ with respect to the Krull topology.
\end{prop}
\begin{proof}
	As $G_K$ maps onto $\Gal(H/K)[2]$, and\footnote{This is independent of the fact that every singular K3 surface has a model over the ring class field.} $G_K \cap \Gal(\bar{K}/H) = \Gal(\bar{K}/H)$, we get the following diagram,
	\[\xymatrix{
		&0 \ar[r] &\Gal(\bar{K}/H) \ar@{=}[d] \ar[r] &G_K \ar@{^{(}->}[d] \ar[r] &\Gal(H/K)[2] \ar@{^{(}->}[d] \ar[r] &0 \\
		&0 \ar[r] &\Gal(\bar{K}/H) \ar[r] &\Gal(\bar{K}/K) \ar[r]^{{\vert_{H}}} &\Gal(H/K) \ar[r] &0
	}
	\]
	from which we extract the short exact sequence
	\[0 \ra G_K \ra \Gal(\bar{K}/K) \ra C(\ko)/C(\ko)[2] \ra 0.\]
	
	The group inclusions $\Gal(\bar{K}/H) \subseteq G_K \subseteq \Gal(\bar{K}/K)$ yield the reversed inclusions of fields $K \subseteq M_K \subseteq H$. Notice that $G_K$ is normal in $\Gal(\bar{K}/K)$. As $G_K = (\vert_H)^{-1} (\Gal(H/K)[2])$, and the restriction map $\vert_H: \Gal(\bar{K}/K) \longra \Gal(H/K)$ is a continuous surjection by Galois theory, $G_K$ is closed and we are done.
\end{proof}

\subsection{The extension $M_K/K$}
We can now use our knowledge of $G_K$ to give a proof of Theorem \ref{thmprim}.

\begin{proof}[Proof of Theorem \ref{thmprim}]
	As $G_K$ is closed and normal in $\Gal(\bar{K}/K)$, we have that
	\[ \Gal(\bar{K}/M_K) = \Gal(\bar{K}/\bar{K}^{G_K})= G_K\]
	and $M_K / K$ is a (finite) Galois extension. The exact sequence
	\[0 \ra G_K \ra \Gal(\bar{K}/K) \ra C(\ko)/C(\ko)[2] \ra 0\]
	tells us that $\Gal(M_K/K) \cong C(\ko)/C(\ko)[2]$, from which we can now cook up the following short exact sequence.
	\[\xymatrix{
		&0 \ar[r] &\Gal(H/K)[2] \ar[r] &\Gal(H/K) \ar[r] &\Gal(M_K/K) \ar[r] &0 \\
		& &\Gal(H/M_K) \ar@{^{(}->}[ru] \ar@{=}[u]^{\text{iso}} & & &
	}
	\]
	
	By genus theory (see \cite{cox13}), there is a short exact sequence
	\[0 \ra C(D)[2] \ra C(D) \ra C(D)^2 \ra 0,\]
	where $C(D)^2$ is the group of squares in the class group $C(D)$ (in fact, it is the principal genus). As $\Gal(H/K) \cong C(D)$, we deduce that
	\[\Gal(M_K/K) \cong C(D)^2,\]
	and in particular that $\#\Gal(M_K/K) = g$, where $g = \#C(D)^2$ is the order of the genus of the transcendental lattice.
\end{proof}

\begin{ex}\label{example}
	Let $D = -23$ and $K = \Q(\sqrt{D})$. The class group of discriminant $D$ is
	\[ C(D) = \Bigg\lbrace \begin{pmatrix} 2 & 1 \\ 1 & 12\end{pmatrix}, \begin{pmatrix} 4 & 1 \\ 1 & 6\end{pmatrix}, \begin{pmatrix} 4 & -1 \\ -1 & 6\end{pmatrix} \Bigg\rbrace.\]
	There is only one genus in $C(D)$ (of order 3), thus we expect a field of moduli of degree 3 over $K$. 
	
	Let $X$ be the singular K3 surface whose transcendental lattice is
	\[P=\begin{pmatrix} 2 & 1 \\ 1 & 12\end{pmatrix}.\]
	A Shioda-Inose structure starting from the self-product of $E$, $E$ being the elliptic curve corresponding to the principal form $P$ in $C(D)$, reveals that $X$ has a model over $\Q(j(P))$. We now show that the field of $K$-moduli is $M_K= K(j(P)) = H(\ko_K)$, which is a degree 3 extension of $K$ by class field theory. Indeed, as $X$ is realized starting from the self-product of $E$, where $E$ corresponds to the principal form $P$, then the transcendental lattice of the conjugate surface by $\sigma \in \Gal(\bar{K}/K)$ is given by 
	\[P^\sigma \ast P^\sigma  = F(\sigma)^{-2},\]
	and this is trivial if and only if $F(\sigma)$ is 2-torsion. However, as $\# \Gal(H/K) =3$, it follows that $F(\sigma)$ is necessarily trivial, and thus $G_K = \ker F = \Gal(\bar{K}/H)$. Therefore, we have proven that $M_K = H$.
	
	We can also look at the K3 surface $Y$ whose transcendental lattice is 
	\[Q:=\begin{pmatrix} 4 & 1 \\ 1 & 6\end{pmatrix}.\]
	By means of a Shioda-Inose structure, $Y$ has a model over $\Q(j_1,j_2)$, where $j_1 := j(P)$ and $j_2:= j(Q)$; notice that
	\[ \Q(j_1,j_2) = K(j_1,j_2) = K(j_2),\]
	as we have considered the Shioda-Mitani model of $Y$ (plus some class field theory considerations). It follows that $H = K(j_2)$, which is a degree 3 extension of $K$, and thus we have that $M_K = H = K(j_2)$. \qed
\end{ex}

\section{Generalization to the imprimitive case}
\subsection{A first look at $G_K$}
We will now treat the case of a singular K3 surface $X$ with imprimitive transcendental lattice $\rT(X) = Q = m Q_0$ ($m>1$). As in the primitive case, we see that it is enough to choose a decomposition of $A$, and to compute the field of moduli in that case. Thus, we now fix a decomposition $A \cong E_1 \times E_2$.

We would like to mimic the techniques used in the primitive case to give an analogous characterization of the field of moduli. The issue at hand is that given a decomposition $A \cong E_1 \times E_2$, the quadratic forms $Q_1$ and $Q_2$ corresponding to $E_1$ and $E_2$ must necessarily lie in class groups with different discriminant by \cite{laface15}. Therefore, we need to use the Dirichlet composition in its generalized sense in order to compute transcendental lattices.

When dealing with decompositions, it is always useful to keep in mind the diagram of orders,
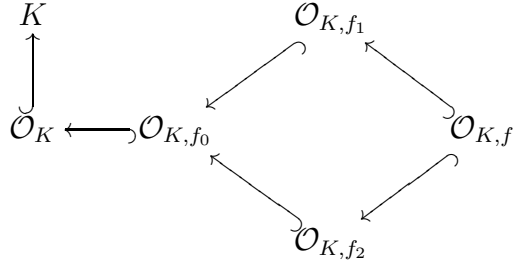
\[\xymatrix{
	K & &\ko_{K,f_1} \ar@{^{(}->}[ld] & \\
	\ko_K \ar@{^{(}->}[u] & \ko_{K,f_0} \ar@{^{(}->}[l] & & \ko_{K,f} \ar@{^{(}->}[lu] \ar@{^{(}->}[ld] \\
	& &\ko_{K,f_2} \ar@{^{(}->}[lu] &
}
\]
and the corresponding one of class groups,
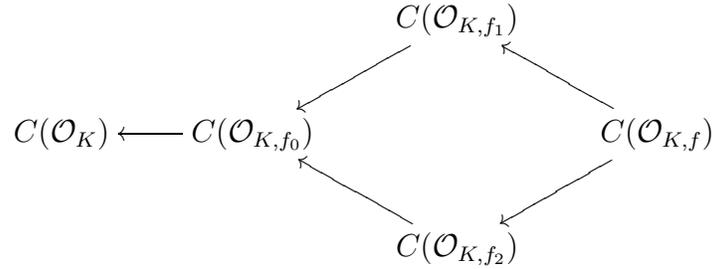
\[\xymatrix{
	& & C(\ko_{K,f_1}) \ar[ld] & \\
	C(\ko_K)  & C(\ko_{K,f_0}) \ar[l] & & C(\ko_{K,f}) \ar[lu] \ar[ld] \\
	& & C(\ko_{K,f_2}) \ar[lu] &
}
\]
where $f_0,f_1,f_2,f$ are such that
\[\lcm(f_1,f_2)=f, \qquad \gcd(f_1,f_2) = f_0, \qquad f^2d_K = \disc \rT(A),\]
and also $[E_1] \in C(\ko_{K,f_1})$ and $[E_2] \in C(\ko_{K,f_2})$. The maps between the above class group are the one induced by extension of ideals; in terms of quadratic forms, these correspond to multiplication by the principal form of the target order: for instance, the reduction map
\[\operatorname{red}_0 : C(\ko_{K,f}) \longra C(\ko_{K,f_0})\]
sends $[Q]$ to $[Q] \circledast [P_0]$, where $\circledast$ is the generalized Dirichlet composition, and $P_0$ is the principal form in $C(\ko_{K,f_0})$. As before, there are maps
\[F_i : \Gal(\bar{K}/K)\longra C(\ko_{K,f_i}) \qquad (i=0,1,2),\]
such that
\[ [Q_i^\sigma] = [F_i(\sigma)]^{-1} \circledast [Q_i] \qquad (i=0,1,2).\]
By use of the generalized Dirichlet composition $\circledast$ and the maps $F_i$ ($i=1,2$), we see that
\begin{align*}
E_1^{\sigma} \times E_2^{\sigma} \cong E_1 \times E_2 \Longleftrightarrow Q_1^{\sigma} \circledast Q_2^{\sigma} = Q_1 \circledast Q_2 \Longleftrightarrow F_1(\sigma) \circledast F_2(\sigma) = P_0.
\end{align*}

The discussion above can be rephrased as follows:

\begin{lemma}\label{lemmaG}
	$G_K = \lbrace \sigma \in \Gal(\bar{K}/K) \, \vert \, F_1(\sigma) \circledast F_2(\sigma) = P_0\rbrace.$
\end{lemma}

In order to go any further, we need to understand the interaction of the maps $F_i$ ($i=0,1,2$). As the class groups are abelian groups, these maps factor through the Galois group of $K^{\rm ab}$, the maximal abelian extension of $K$. We get maps (again called $F_i$ by abuse of notation)
\[F_i : \Gal(K^{\rm ab}/K) \longra C(\ko_{K,f_i}).\]
Here is where the theory of id\'{e}les comes into play, picturing the behaviour of these maps in their totality.

\subsection{Compatibility condition for the maps $F_i$}
The idea is inspired by a paper of Sch\"utt \cite{schuett07}: given a singular abelian surface $A$, among all decompositions that we can choose, there is one that behaves better that the others, namely the decomposition that Shioda and Mitani used to prove the surjectivity of the period map for singular abelian surfaces \cite{shioda-mitani74}. 

To the reader's convenience, we briefly recall this construction. Letting $A$ be a singular abelian surface of transcendental lattice $\rT(A) \cong (a,b,c)$, Shioda and Mitani showed that $A \cong E_{\tau} \times E_{a\tau+b}$, where
\[\tau := \tau (Q) = \frac{-b+\sqrt{D}}{2a}.\]
In particular, $E_{a\tau+b}$ always corresponds to the principal form in the class group of discriminant $D = \disc \rT(A)$, and $E_\tau$ instead corresponds to the quadratic form $\rT(A)_0$, the primitive part of $\rT(A)$.

Let us assume $A \cong E_1 \times E_2$ is the Shioda-Mitani decomposition: if $\rT(A) = Q = mQ_0$, then $E_1$ corresponds to the quadratic form $Q_0 \in C(D_0)$ and  $E_2$ corresponds to the principal form $P \in C(D)$.  Notice that we also have $A \cong \C / \ga \times \C/ \ko_{K,f}$, for $\ga \in C(\ko_{K,f_0})$, and thus the proof of \cite[Theorem 5.4]{schuett07} shows in particular that, for $\sigma \in \Gal(\bar{K}/K)$
\[A^\sigma \cong E_1^\sigma \times E_2^\sigma \cong \C / s^{-1}\ga \times \C/s^{-1}\ko_{K,f} \cong \C / s^{-2}\ga \times \C/\ko_{K,f},\]
where, as $s$ varies in $\II_K$, $(s^{-1}\ga)^2 = (s^{-1}Q_0)^2$ spans the whole genus of $Q_0$ in $C(D_0)$. This ultimately suggests that we look at elements of the form $s^{-1}\ko$, as their squares span the principal genus of a class group, and characterize the transcendental lattice as it moves in its genus.

To do so, suppose we are given an order $\ko \subset K$, the map
\[F: \Gal(\bar{K}/K) \longra C(\ko)\]
factorizes through a map
\[F :\Gal(K^{\rm ab}/K) \longra C(\ko).\]
For $\sigma \in \Gal(K^{\rm ab}/K)$, $F(\sigma)$ has the property $[E^\sigma] = [\C / F(\sigma)^{-1} \cdot \ga]$ independently of the chosen $E = \C/\ga \in \ke ll(\ko)$. By the Main Theorem of CM, there exists an id\'{e}le $s \in \II_K$ such that $\phi_K(s) = \sigma$ and 
\[ [E^\sigma ] = [\C / F(\sigma)^{-1} \cdot \ga] = [\C / s^{-1} \ga].\]
As $s^{-1}\ga = (s^{-1}\ko)\cdot \ga$, we can identify $[s \ko]  = [F(\sigma)]$. Now let $\ko_0$ be another order in $K$, $\ko \subset \ko_0 \subset K$, and consider the following diagram.

\[\xymatrix{
	& & & & &C(\ko)\ar[dd]^{\operatorname{red_0}}\\
	&\II_K \ar[rr]^{\phi_{K}} & &\Gal(K^{\rm ab}/K) \ar[rru]^F \ar[rrd]_{F_0} & & \\
	& & & & &C(\ko_0)
}
\]

We would like to show that the triangle on the right-hand side is indeed commutative. For $\sigma \in K^{\rm ab}$, we have the identifications
\[ [F(\sigma)] = [s \ko] \qquad \text{and} \qquad [F_0(\sigma)] = [s\ko_0],\]
which are a consequence of the Main Theorem of CM. Notice that this uses the fact that the Main Theorem of CM holds for all elliptic curves with CM in any order in $K$ at once. By looking at every rational prime $p$, one checks that $(s\Lambda)\cdot \Lambda' = s(\Lambda \cdot \Lambda ')$, for two lattices $\Lambda$ and $\Lambda'$ in $K$ (see \cite{shimura71}). In particular, after noticing that $\Lambda$ and $s \Lambda$ have the same endomorphism ring, we get $[s\ko]\circledast[\ko_0] = [s\ko_0]$. We have proven the following compatibily condition
\begin{lemma}\label{lemmacomm}
	Under the assumptions above,
	\[[F_0(\sigma)] = [F(\sigma)] \circledast [P_0],\]
	or equivalently $\operatorname{red}_0 \circ F = F_0$.
\end{lemma}

This proves the commutativity of the triangle in the diagram above, and thus we are now ready to prove a characterization theorem for the field of $K$-moduli also in the imprimitive case. 

\subsection{Completion of the proof}
In Lemma \ref{lemmaG}, we showed that
\[E_1^{\sigma} \times E_2^{\sigma} \cong E_1 \times E_2 \Longleftrightarrow  F_1(\sigma) \circledast F_2(\sigma) = P_0.\] 
Now, $F_1(\sigma) \circledast F_2(\sigma)$ lives in $C(\ko_{K,f_0})$ so we can multiply by the principal form $P_0$, and, by commutativity and Lemma \ref{lemmacomm}, the last condition above is equivalent to $F_0(\sigma)^2 = P_0$, i.e.~$F_0(\sigma) \in C(\ko_{K,f_0})[2]$. Therefore we get, in analogy to the primitive case:

\begin{prop}
	$G_K = F_0^{-1}(C(\ko_{K,f_0})[2])$.
\end{prop}

Now, the same argument used in the primitive case (replacing every occurrence of $H$ with $H_0$, the ring class field of $\ko_{K,f_0}$), yields the following result, which extends Theorem \ref{thmprim} to the imprimitive case.
\begin{thm}\label{relmoduli}
	Let $X$ be a singular K3 surface with transcendental lattice $T(X)= Q = mQ_0$, and let $H_0$ be the ring class field of $\ko_{K,f_0}$, the order of discriminant $\disc Q_0$. Then the field of $K$-moduli is
	\[M_K = \bar{K}^{G_K}, \qquad G_K = (\vert_{H_0})^{-1}\Gal(H_0/K)[2];\]
	it is a Galois extension of $K$ of degree $g$, $g$ being the order of the genus of the transcendental lattice of $X$.
\end{thm}

\section{The absolute field of moduli}\label{abs_field}
So far, we have studied the field of $K$-moduli of a singular K3 surface $X$, $K$ being the CM field of $X$. Now, we want to move our attention to the \textit{absolute field of moduli} $M_\Q$, by which we mean the field of $\Q$-moduli. We will proceed as in the case of $M_K$.

Let us recall that the absolute field of moduli of $X$ is the field $M_\Q:= \C^{G_\Q}$, where 
\[G_\Q = \lbrace \sigma \in \Gal( \C/\Q) \, \vert \, X^\sigma \in [X] \rbrace.\]
The proof of Lemma \ref{FM} shows that we can equivalently define the field of moduli $M_\Q$ to be the subfield of $\bar{\Q}$ which is fixed by the group
\[G_\Q = \lbrace \sigma \in \Gal( \bar{\Q}/\Q) \, \vert \, X^\sigma \in [X] \rbrace.\]
As $G_K$ is the subgroup of elements of $G_\Q$ whose restriction to $K$ is trivial, we have the following commutative diagram,
\[\xymatrix{
	&0 \ar[r] &G_K \ar@{^{(}->}[d] \ar[r] &G_\Q \ar@{^{(}->}[d] \ar[r] &C \ar[d] \ar[r] &0 \\
	&0 \ar[r] &\Gal(\bar{\Q}/K) \ar[r] &\Gal(\bar{\Q}/\Q) \ar[r]^{{\vert_{K}}} &\Gal(K/\Q) \ar[r] &0
}
\]
where $C$ is simply the quotient group $G_\Q / G_K$ (notice that $G_K$ is normal in $G_\Q$). We have the following:

\begin{prop}\label{prop_6.1} $C \cong \Gal(K / \Q)$.
\end{prop}

\begin{proof}
This is equivalent to $G_K$ being an index-two subgroup of $G_\Q$. In fact, it is enough to show that $G_{\Q} \setminus G_K \neq \emptyset$. Indeed, assume there exists an element $\sigma \in G_\Q \setminus G_K$, thus $X^\sigma \cong X$ and $\sigma \vert _K \neq \id_K$. Also, notice that $\sigma^{-1} \in G_\Q \setminus G_K$. If $\tau \in G_\Q \setminus G_K$ is another such element, then $\sigma \tau^{-1} \in G_K$, which means $\bar{\sigma} = \bar{\tau} \in C$. This implies that $C \cong \Z / 2\Z$, and thus we can identify the quotient group $C$ with $\Gal(K/\Q)$.
	
	We now prove that $G_{\Q} \setminus G_K \neq \emptyset$. If $\rT(X)$ is 2-torsion, then the complex conjugation automorphism $\iota$ is an element in $G_\Q \setminus G_K$. Therefore we can assume $\rT(X)$ is not 2-torsion. In this case, we will build an element $\alpha \in G_\Q \setminus G_K$. Suppose $X$ is given as the Shioda-Inose K3 surface of a product $E_1 \times E_2$, and let $Q_i \in C(\ko_{K,f_i})$ be the corresponding quadratic forms (for $i=1,2$). Then $\rT(X) = m Q_0$ with $Q_0 = Q_1 \circledast Q_2 \in C(\ko_{K,f_0})$, $m \geq 1$ and $f_0 = \gcd(f_1,f_2)$. Then, for $\sigma \in \Gal(\bar{\Q}/K)$, consider $X^{\iota \sigma} = (X^\sigma)^\iota$: by looking at the transcendental lattice (as a quadratic form), we see that
	\begin{align*}
		\rT(X^{\iota \sigma}) &= \rT(X^\sigma)^{-1} = \big[ m(F_0(\sigma)^{-2} \circledast Q_1 \circledast Q_2) \big]^{-1} = m \big[ F_0(\sigma)^2 \circledast (Q_1 \circledast Q_2)^{-1} \big] .
	\end{align*}
	Notice that we have used the fact that complex conjugation is insensitive to the index of primitivity of $\rT(X)$ and acts on it by inversion (i.e.~by multiplication of the non-diagonal entries of $\rT(X)$ by $-1$, after choosing a basis). At this point, we choose $\sigma \in \Gal(\bar{\Q}/K)$ such that $F_0(\sigma) = Q_0$. By Lemma \ref{lemmacomm}, this choice is compatible with the reduction maps $\operatorname{red}_0$ of quadratic forms. Such $\sigma$ yields an element $\alpha:= \iota \sigma$ as required.

\end{proof}

As a consequence, $M_K / M_\Q$ is a Galois extension with group $C \cong \Gal(K/\Q)$, and thus $M_K \supsetneqq M_\Q$. By multiplicativity of degree, the extensions $M_K/ K$ and $M_\Q / \Q$ have the same degree. 

We remark that $K$ is not contained in $M_\Q$, and thus $M_\Q \cap K = \Q$. Indeed, if this were the case, then we would have $G_\Q \subset \Gal(\bar{\Q}/K)$. If $\rT(X)$ is 2-torsion, then $\iota \in G_\Q \setminus \Gal(\bar{\Q}/K)$, and so we get a contradiction. If $\rT(X)$ is not 2-torsion, then any $\tau \in G_\Q \setminus G_K$ yields the same contradiction.

As a consequence of the discussion above, we have the following result.
\begin{thm}\label{absmoduli}
	Let $X$ be a singular K3 surface. Its absolute field of moduli $M_\Q$ is an index-two subfield of the field of $K$-moduli $M_K$. Moreover, $M_\Q$ is an extension of $\Q$ of degree
	\[ [M_\Q : \Q] = [M_K : K] = g,\]
	$g$ being the genus of $T(X)$. In general, it is not a Galois extension of $\Q$.
\end{thm}

\begin{figure}
	\[
	\xymatrix{
		&\bar{\Q} & & & \\
		& &H\ar@{-}[lu] & & \\
		& & &M_K \ar@{-}[lu] & \\
		&K \ar@{-}[uuu] \ar@{-}[ruu] \ar@{-}[rru] & & &M_\Q \ar@{-}[lu] \\
		& &\Q \ar@{-}[lu] \ar@{-}[rru] & &
	}
	\]
	\caption{Relative and absolute field of moduli}\label{ext}
\end{figure}
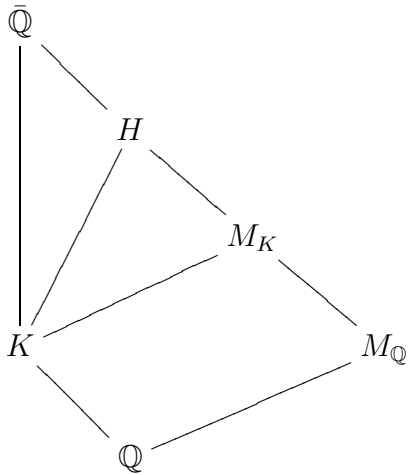

In particular, we have the diagram of field extensions in Figure \ref{ext}, where all extensions are Galois, except possibly for $M_\Q / \Q$.

\begin{rmk}
	We would like to point out that Figure \ref{ext} recovers the picture of the case of elliptic curves with CM in $K$. For an elliptic curve $E$ with CM in an imaginary quadratic field $K$, $M_\Q = \Q(j(E))$ and $M_K = K(j(E)) = H$, $H$ being the ring class field of the elliptic curve $E$ (we are implicitly using the fact that elliptic curves correspond to quadratic forms). The equality $M_K = H$ is explained by the fact that the field of $K$-moduli coincides with the minimal field of definition, for every elliptic curve with CM in $K$. 
\end{rmk}

The final statement of Theorem \ref{absmoduli} is that the extension $M_\Q / \Q$ is not Galois in general; the following example shows an occurrence of this phenomenon.

\begin{ex}[Example \ref{example} reloaded]\label{counterex}
	We compute the absolute field of moduli for the K3 surface $X$. Our results tell us that $M_\Q$ must be an extension of degree 3 of $\Q$. In this case, $X$ has a model over $\Q(j(P))$, which is a degree 3 extension of $\Q$, so it follows that the absolute field of moduli $M_\Q$ is indeed $\Q(j(P))$ itself, which agrees with Theorem \ref{absmoduli}. Now, the class polynomial $H_{\ko_K}(T)$ of the order $\ko_K$ has $j(P)$, $j(Q)$ and $j(Q^{-1})$ as roots; $j(P)$ is real, while $j(Q)= \overline{j(Q^{-1})}$. It follows that the extension $\Q(j(P))/\Q$ cannot be Galois.
	
	For a more interesting example, we look at the K3 surface $Y$, and we recall that
	\[\Q(j_1,j_2) = K(j_1,j_2) = K(j_2),\qquad [\Q(j_1,j_2):\Q]=6,\]
	if we consider a Shioda-Inose model for $Y$. As the field of moduli $M_\Q$ in contained in every field of definition for $Y$, and must have degree 3 by Theorem \ref{absmoduli}, we must find an element $\alpha \in \Gal(H/\Q)$ which leaves the modulus invariant. If $\iota$ denotes the (restriction of the) complex conjugation automorphism and $\sigma \in \Gal(\bar{K}/K)$ is an element such that $F(\sigma) = Q^{-1}$, then $\alpha := \sigma \iota$ satisfies this condition. Therefore $M_\Q$ is the subfield of $H$ which is fixed by the group generated by $\alpha$: this group has order 2, and thus we get that $M_\Q$ is an extension of $\Q$ of degree 3, as expected.\qed
\end{ex}

\section{Further questions}

\subsection{Non-finiteness of singular K3 surfaces}
The following discussion is inspired by the following striking result of \v{S}afarevi\v{c} on the finiteness of singular K3 surfaces with bounded field of definition:

\begin{thm}[Theorem 1 of \cite{shafarevich96}]
	Let $n$ be a positive integer. There exist finitely many singular K3 surfaces with a model over a number field $K$ of degree $[K:\Q] \leq n$.
\end{thm}

This result says that we can use the degree of the field of definition to stratify $\Sigma^{\rm K3}$, and that each stratum contains finitely many elements only: the $n$\textsuperscript{th} stratum is defined as
\[ \Sigma^{\rm K3} (n) := \lbrace [X] \in \Sigma^{K3} \, : \, \text{$X$ has a model over $K$, $[K:\Q] \leq n$} \rbrace.\]

One might wonder whether a similar result holds for the field of moduli in place of the field of definition. We will now see that this is not the case.

\begin{prop}\label{red}
	Let $X$ and $Y$ be two singular K3 surfaces such that $\rT(X)$ is primitive (as a quadratic form) and $\rT(Y) = n \rT(X)$, for some $n \in \IN$. Then, $X$ and $Y$ have the same field of $K$-moduli, $K$ being the CM field of $X$ and $Y$.
\end{prop}

\begin{proof}
	The argument used in proving Theorem \ref{relmoduli} shows, in particular, that the ring class field $H_0$ only depends on the discriminant of the primitive part of the transcendental lattice. In the situation at hand, $X$ and $Y$ would both lead to the same ring class field, and the result is then a consequence of Theorem \ref{relmoduli}.
\end{proof}

\begin{prop}\label{red1}
	Let $X$ and $Y$ be two singular K3 surfaces whose transcendental lattices are primitive and lie in the  same class group (as quadratic forms). Then, $X$ and $Y$ have the same field of $K$-moduli, $K$ being the CM field of $X$ and $Y$.
\end{prop}

\begin{proof}
	Same as for Proposition \ref{red}.
\end{proof}

As a corollary, we get that

\begin{cor}\label{cor_abs}
	Let $X$ and $Y$ be two singular K3 surfaces such that the primitive parts of $\rT(X)$ and $\rT(Y)$ lie in the same class group (as quadratic forms). Then, $X$ and $Y$ have the same field of $K$-moduli, $K$ being the CM field of $X$ and $Y$.
\end{cor}

In particular, this shows that bounding the degree of the (relative) field of moduli is not enough to have a stratification of $\Sigma^{\rm K3}$ in strata containing finitely many elements only. In fact, we have shown that for each possible field of $K$-moduli, there exist infinitely many singular K3 surfaces with that field of $K$-moduli. This non-finiteness result holds true also if we replace the relative field of moduli with the absolute one: in fact, it is enough to fix a primitive quadratic form $Q$ such that $h(\disc Q) =1$; then 
\[ \# \lbrace [X] \in \Sigma^{\rm K3} \, : \; \rT(X) = m Q, \ m \in \IN \rbrace = + \infty,\]
and all K3 surfaces in the set above have clearly $\Q$ as absolute field of moduli.

\subsection{Explicit fields of $K$-moduli}
We can still ask questions such as: which fields can appear as the field of $K$-moduli of a singular K3 surface? To answer such a question, Theorem \ref{relmoduli} and its description of the field of moduli does not help us. The ideal situation would be to describe $M_K$ as the subfield of a finite extension of $K$ fixed by a (finite) group. This would also allow us to explicitly describe this field with the aid of a computer algebra system.

In consequence of Proposition \ref{red}, we can restrict ourselves to working with singular K3 surfaces whose transcendental lattice is primitive as a quadratic form; thus, let $X$ be such a singular K3 surface. Consider $X$ as obtained by a singular abelian surface $A$ in its Shioda-Mitani model $A \cong E_Q \times E_P$. Let us remind the reader that $Q$ and $P$ belong to the same class group, exactly because the transcendental lattice is primitive. Then, a result of Sch\"utt \cite{schuett07} implies that $X$ has a model\footnote{More generally it has a model over the field $\Q(j_1,j_2)$, which does not always coincide with the ring class field $H$.} over the ring class field
\[H:=K(j_1,j_2) = K(j_2), \qquad j_k:=j(E_k) \quad (k=1,2).\]
This model is particularly nice as the extension $H/K$ is Galois by class field theory. It is clear that $\Gal(\bar{K}/H) \subseteq G_K$, because of the existence of a model over $H$. Also, the arguments in Section \ref{prim} yield a proof of the following result:

\begin{prop}\label{explicit}
	$M_K = H^{\Gal(H/K)[2]}.$
\end{prop}

\begin{proof}
	Let us consider the restriction map $\vert_H : \Gal(\bar{K}/K) \longra \Gal(H/K)$. 	By the existence of a model over $H$, $G_K$ maps onto the following subgroup of $\Gal(H/K)$: 
	\[G_K\vert_H := \lbrace \sigma \in \Gal(H/K) \, \vert \, X^\sigma \cong X \rbrace.\]
	The proof of Proposition \ref{fundcond} shows in particular that $G_K \vert_H = \Gal(H/K)[2]$. Now the proof follows after a direct check. 
\end{proof}

This last result allows us to explicitly compute the field of moduli of a given singular K3 surface. We can have a computer algebra system run this sort of computations for us, but in order to do so, we have to reduce to isolate a finite number of cases at the time. To this end, Proposition \ref{red} enables us to project $\Sigma^{\rm K3}$ onto
\[ \Sigma^{\rm K3}_{\text{prim}} := \lbrace [X] \in \Sigma^{\rm K3} \, : \, \rT(X) \text{ is primitive} \rbrace,\]
by forgetting the index of primitivity of the transcendental lattice. Analogously to the situation of \cite{shioda-mitani74}, there is a 1:1 correspondence
\[ \Sigma^{\rm K3}_\text{prim} \longleftrightarrow \kq_0^+ / \SL_2(\Z),\]
where $\kq_0^+$ is the subset of $\kq^+$ containing primitive quadratic forms only. Class group theory implies that
\[ \kq_0^+ / \SL_2(\Z) \cong \bigsqcup_{\substack{\text{$K$ quadratic imaginary field} \\ f \in \IN}} C(\ko_{K,f}),\]
and thus we can bound $\Sigma^{\rm K3}_\text{prim}$ by bounding the orders in the quadratic imaginary fields. This can be achieved, for instance, by bounding the discriminant or the class number. Such constraint gives a stratification of $\Sigma^{\rm K3}_\text{prim}$ whose strata contain finitely many elements only, and we can therefore run the computations in a finite, perhaps long, time.

\subsection{Invariance of the field of $\Q$-moduli}
Let $X$ and $X$ be singular K3 surfaces such that their transcendental lattices have the same primitive part, i.e.~$\rT(X) = mQ_0$, $\rT(Y)=m'Q_0'$ and $Q_0 = Q_0'$ (primitive). From previous discussions, we know that $X$ and $X'$ will have the same field of $K$-moduli, with $K = \Q(\sqrt{\disc \rT(X)})$. We would like to study the analogous question for the field of $\Q$-moduli, that is whether the field of $\Q$-moduli is independent of the index of primitivity.

\begin{prop}
	Let $X$ and $Y$ be two singular K3 surfaces such that $\rT(X)$ is primitive (as a quadratic form) and $\rT(Y) = m \rT(X)$, for some $m \in \IN$. Then, $X$ and $Y$ have the same field of $\Q$-moduli, i.e.~$M_\Q(X) = M_\Q(Y)$. 
\end{prop}

\begin{proof}
The discussion in Section \ref{abs_field} shows that for  $G_\Q(X) = \langle G_K(X),\alpha_X \rangle$ (as a group), where $\alpha_X$ is a suitable element in $\Gal(\bar{\Q}/K)$. We will now show that $\alpha_X = \alpha_Y$, where $\alpha_Y$ is the analogous element of $\Gal(\bar{\Q}/K)$ for $Y$. 

Indeed, if we write $\alpha_X = \iota \circ \sigma_X$, with $\sigma_X \in \Gal(\bar{\Q}/K)$ as in the proof of Proposition \ref{prop_6.1}\footnote{When $\rT(X)$ is 2-torsion (and thus also $\rT(Y)$ is), $\alpha_X = \alpha_Y = \iota$, $\iota$ being the complex conjugation automorphism, and thus $\sigma_X = \text{id}_X$.}, it is straightforward to see that
\[ \rT(Y^{\alpha_X}) = \rT(Y^{\alpha_X})^{-1} = m\big[F_0(\sigma_X)^{-2} \circledast \rT(X) \big]^{-1} = m \big[ F_0(\sigma_X)^{2} \circledast \rT(X)^{-1} \big] = m \rT(X) = \rT(Y), \]
where we have used that $\rT(X)$ is primitive and that $F_0(\sigma_X) = \rT(X)$ (as a quadratic form -- see proof of Proposition \ref{prop_6.1}). 
\end{proof}

This result is clearly the analogue of Proposition \ref{red} for the field of $\Q$-moduli. Since conjugation by an automorphism acts as taking the conjugate field on the field of moduli (see Section \ref{section_3}), we cannot have such an analogue for Proposition \ref{red1}. However, the following analogue of Corollary \ref{cor_abs} holds.

\begin{cor}
	Let $X$ and $Y$ be two singular K3 surfaces. Assume that there exists $\sigma \in \Gal(\bar{\Q}/\Q)$ such that $\rT(X^\sigma)_0 = \rT(Y)_0$ as quadratic forms, $\rT(Y)_0$ being the primitive part of $\rT(Y)$ (and similarly for $X^\sigma$). Then $M_\Q(Y) = \sigma \big( M_\Q(X) \big)$.
\end{cor}

\begin{proof}
	Use the fact that $M_\Q(X^\sigma) = \sigma \big( M_\Q(X) \big)$ together with
	\[ \rT(X^\sigma) = m \rT(X^\sigma)_0 = m \big[ F_0(\sigma)^{-2} \circledast \rT(X)_0 \big].\]
\end{proof}

Finally, thanks to Proposition \ref{explicit}, it is straightforward to describe $M_\Q$ explicitly.

\begin{prop}
	$M_\Q = H^{\Gal(H/K)[2] \times \Gal(K/\Q)}$.
\end{prop}

\bibliographystyle{plain}
\bibliography{bib}{}
\end{document}

%% file: pack.tex
\usepackage{amsmath}
\usepackage{amsfonts}
\usepackage{amssymb}
\usepackage{amsthm}
\usepackage{mathtools}
\usepackage{mathrsfs}
\usepackage{float}
\usepackage[all,cmtip]{xy}
      \theoremstyle{definition}
\newtheorem{defi}{Definition}[section]
\theoremstyle{plain}
\newtheorem{thm}[defi]{Theorem}

\newtheorem{prop}[defi]{Proposition}
\newtheorem{cor}[defi]{Corollary}
\newtheorem{lemma}[defi]{Lemma}
\theoremstyle{remark}

\newtheorem{rmk}[defi]{Remark}
\newtheorem{ex}[defi]{Example}

\theoremstyle{definition}

\newcommand{\ra}{\rightarrow}
\newcommand{\longlra}{\longleftrightarrow}

\newcommand{\longra}{\longrightarrow}

\newcommand{\ke}{{\mathcal E}}

\newcommand{\ko}{{\mathcal O}}

\newcommand{\kq}{{\mathcal Q}}

\newcommand{\C}{{\mathbb C}}

\newcommand{\II}{{\mathbb I}}

\newcommand{\IN}{{\mathbb N}}

\newcommand{\IP}{{\mathbb P}}
\newcommand{\Q}{{\mathbb Q}}

\newcommand{\Z}{{\mathbb Z}}

\newcommand{\ga}{{\mathfrak a}}
\newcommand{\gb}{{\mathfrak b}}

\newcommand{\gm}{{\mathfrak m}}

\newcommand{\rH}{{\rm H}}
\newcommand{\rT}{{\rm T}}

\newcommand{\Aut}{\operatorname{Aut}}

\newcommand{\GL}{\operatorname{GL}}
\newcommand{\SL}{\operatorname{SL}}

\newcommand{\Gal}{\operatorname{Gal}}

\newcommand{\NS}{\operatorname{NS}}

\newcommand{\id}{\operatorname{id}}

\newcommand{\ord}{\operatorname{ord}}
\newcommand{\Km}{\operatorname{Km}}
\newcommand{\Div}{\operatorname{Div}}

\newcommand{\rank}{\operatorname{rank}}   
\newcommand{\disc}{\operatorname{disc}}

\newcommand{\lcm}{\operatorname{lcm}}

  \makeatletter
\newcommand{\xdashrightarrow}[2][]{\ext@arrow 0359\rightarrowfill@@{#1}{#2}}
\newcommand{\xdashleftarrow}[2][]{\ext@arrow 3095\leftarrowfill@@{#1}{#2}}
\newcommand{\xdashleftrightarrow}[2][]{\ext@arrow 3359\leftrightarrowfill@@{#1}{#2}}
\def\rightarrowfill@@{\arrowfill@@\relax\relbar\rightarrow}
\def\leftarrowfill@@{\arrowfill@@\leftarrow\relbar\relax}
\def\leftrightarrowfill@@{\arrowfill@@\leftarrow\relbar\rightarrow}
\def\arrowfill@@#1#2#3#4{%
  $\m@th\thickmuskip0mu\medmuskip\thickmuskip\thinmuskip\thickmuskip
   \relax#4#1
   \xleaders\hbox{$#4#2$}\hfill
   #3$%
}
\makeatother
